\documentclass{amsart}

\usepackage{amsfonts,amsmath,amssymb,amsthm,latexsym,enumerate,mathrsfs}
\usepackage{braket}
\usepackage{mathtools}
\usepackage{tikz-cd}
\usepackage{hyperref}
\usepackage{enumitem}
\usepackage{cleveref}
\usepackage{setspace}
\usepackage{xcolor}
\usepackage{textcomp}
\usepackage{comment}

\usepackage{blindtext}
\usepackage[a4paper]{geometry}
\newcommand{\Z}{\mathbb{Z}}
\newcommand{\R}{\mathbb{R}}
\newcommand{\C}{\mathbb{C}}

\newcommand{\Hess}{\mathrm{Hess}}
\renewcommand{\d}{\nabla}
\newcommand{\w}{\wedge}
\newcommand{\pp}{\partial}
\newcommand{\il}{\langle}
\newcommand{\ir}{\rangle}

\newtheorem{thm}{Theorem}[section]
\newtheorem{lmm}[thm]{Lemma}
\newtheorem{prop}[thm]{Proposition}
\newtheorem{rmk}[thm]{Remark}
\newtheorem{cor}[thm]{Corollary}

\newtheorem{ex}[thm]{Example}
\newtheorem{remark}[thm]{Remark}

\begin{document}

\setstretch{1.05}

\title[Short Title]{Global Hypersurfaces of Section for Geodesic Flows on Convex Hypersurfaces}

\author{Sunghae Cho and Dongho Lee}

\begin{abstract}
    We construct a global hypersurface of section for the geodesic flow of a convex hypersurface in Euclidean space admits an isometric involution.
    This generalizes the Birkhoff annulus to higher dimensions.
\end{abstract}

\maketitle

\section{Introduction}

Global surfaces of section, an idea introduced by Poincar\'{e} in his work on celestial mechanics \cite{Poincare_87} and also explored by Birkhoff \cite{Birkhoff_66}, feature prominently in the literature on the 3-dimensional dynamics.
They allow us to reduce the dynamics of vector fields on 3-manifolds to the dynamics of surface diffeomorphisms.
Ghys \cite{Ghys_09} called the existence of a global surface of section as a paradise for dynamicists, since it eliminates technical difficulties and allows to investigate the pure nature of dynamical systems.

Beyond their versatile utility, the existence of the global surfaces of section already gives us some information about the dynamics.
Unless the given 3-manifold fibers over a circle, a global surface of section must have a boundary which is a set of periodic orbits.
Hence the existence of periodic orbits is an essential obstruction, which cannot be overcome easily.
Kuperberg \cite{Kuperberg_94} introduced an example of a nonvanishing vector field on $S^3$ without periodic orbits.
Ginzburg \cite{Ginzburg_99} has given examples of the Hamiltonian system without periodic orbits.
The horocycle flow on unit cotangent bundle of higher genus surfaces provides another class of Hamiltonian flows without periodic orbits.

This obstruction can be removed in the case of Reeb dynamics in dimension three.
The \emph{Weinstein conjecture} \cite{Weinstein_79} asserts the existence of at least one periodic Reeb orbit on any compact contact manifold.
Hofer \cite{Hofer_93} proved the Weinstein conjecture on $S^3$, and Taubes \cite{Taubes_07} proved the Weinstein conjecture in dimension three in a full generality using embedded contact homology.
The work of Hofer was followed by the joint work with Wysocki and Zehnder \cite{Hofer_Wysocki_Zehnder_98} constructing disk-like global surface of section on a dynamically convex $S^3$.
Recent developments along this way can be found in \cite{Hryniewicz_Salomao_11}, \cite{Hryniewicz_Momin_Salomao_15}, and \cite{Hryniewicz_Salomao_18}.
\cite{Colin_Dehornoy_Hryniewicz_Rechtman_23}, \cite{Contreras_Mazzucchelli_22} and \cite{Contreras_Knieper_Mazzucchelli_Schulz_22} contain results about the existence of a Birkhoff section for Reeb flows in dimension three, a generalization of a global surface of section that allows immersed boundary.

The notion of global surface of section can be generalized to higher dimensions as mentioned by Birkhoff \cite{Birkhoff_66}.
However, there are only few results known in higher dimensional dynamics even in the case of Reeb flows, due to the difficulty caused by the instability of the boundary of global hypersurface of section.
Moreno and van Koert \cite{Moreno_van_Koert_22a} resolved this problem by imposing a symmetry which guarantees the existence of high dimensional invariant set, and provided a concrete example of a global hypersurface of section arising from the restricted three-body problem.

Our main result concerns the existence of global hypersurfaces of section of the geodesic flows on compact convex hypersurfaces in Euclidean space which satisfying certain conditions.

\begin{thm}
    Let $M\subset\R^{n+1}$ be a regular closed hypersurface and $N\subset M$ be a codimension 1 submanifold.
    Assume that there exists a tubular neighborhood $\nu(N)$ of $N$ and an isometric involution $i:\nu(N)\to\nu(N)$ whose fixed point locus is $N$.
    In addition, assume that $M$ has positive sectional curvature.
    Then the geodesic flow on $ST^*M$ admits a global hypersurface of section
\[
P = \left\{(x,y)\in ST^*M : x\in N, \langle y,\nu_x\rangle \geq 0 \right\},
\]
where $\nu$ is a fixed normal vector field of $N$ with respect to $M$. Moreover, the return map $\Psi:\mathring{P}\to\mathring{P}$ extends smoothly to the boundary of $P$.
\end{thm}

We note that the positivity of sectional curvature is equivalent to the convexity of the hypersurface.
We first prove the existence part in \Cref{Open_Book_Existence_2}, and
\Cref{Smooth_Extension_to_Boundary} establishes that we can extend the return map smoothly to the boundary.
Note that if the dimension of $M$ is 2, our theorem reduces to a special case of the theorem of Birkhoff \cite{Birkhoff_66} which established the existence of a global surface of section of geodesic flow of a 2-sphere with positive curvature.
In this sense, our result can be regarded as a generalization of the Birkhoff annulus to higher dimensions.

\subsection{Acknowledgement}
Sunghae Cho and Dongho Lee were supported by National Science Foundation of Korea Grant NRF2023005562 funded by the Korean Government.

\section{Preliminaries}

\subsection{Global Hypersurfaces of Section and Open Book Decompositions}

Let $Y$ be a smooth closed manifold, i.e.\,\,a compact manifold without boundary, and $X$ be a non-vanishing vector field on $Y$. A \textbf{global hypersurface of section} for $X$ is an embedded submanifold $P \subset Y$ of codimension 1 with (possibly empty) boundary $\partial P = B$ such that
\begin{itemize}[label={$\cdot$}]
    \item the vector field $X$ is transverse to the interior $\mathring{P}$ of $P$,
    \item the boundary $B$ is $X$-invariant;
    in other words, $X$ is tangent to $B$,
    \item for any $p \in Y$, there exist $t_+, t_- > 0$ such that $Fl^X_{t_{+}}(p), Fl^X_{-t_-}(p) \in P$.
\end{itemize}
If $P$ is a global hypersurface of section, we can define the \textbf{(first) return time} $\tau_p$ for each $p \in \mathring{P}$ by $\tau_p = \min\{t>0 : Fl^X_t(p) \in P\}$ and the \textbf{(first) return map} by $\Psi(p) = Fl^X_{\tau_p}(p)$.

\begin{ex}\label{Standard_Sphere_Example}
\rm
Let $Y$ be $ST^*S^n$, the unit cotangent bundle of the $n$-sphere with the standard metric, and $X$ be the unit geodesic vector field on $Y$. Consider $Y$ as a subset of $T^*\R^{n+1} \simeq \R^{n+1} \times \R^{n+1}$ and use the coordinates $(x,y) = (x_0,\cdots,x_n,y_0,\cdots,y_n)$. Consider the submanifold
$$
P = \left\{(x,y) \in Y : x_0=0, y_0 \geq 0\right\},
$$
which is the set of upward directions on the equator.
Then one can see that $P$ is a global hypersurface of section, the boundary of $P$ is the unit cotangent bundle of the equator of $S^n$, and the first return map is the identity.
\end{ex}

Open book decompositions are closely related to the global hypersurfaces of section.
An \textbf{open book decomposition} on a closed manifold $Y$ is a pair $(B,\pi)$ of a codimension 2 closed submanifold $B$ and a map $\pi$ from $Y\setminus B$ to $S^1\subset\C$ which satisfies the following.
\begin{itemize}[label={$\cdot$}]
    \item The normal bundle of $B$ is trivial. We call $B$ the \textbf{binding}.
    Let $\xi:B\times D^2\to \nu(B)$ be a fixed trivialization of the normal bundle, which is identified with a tubular beighborhood, of $B$.
    \item The map $\pi$ is a fiber bundle such that $(\pi\circ\xi)(b;r,\theta)=e^{i\theta}$ on $\nu(B)\setminus B$, where $(r,\theta)$ is a polar coordinate on $D^2$.
    We call the closure of each fiber $\overline{\pi^{-1}(e^{i\theta})}=P_\theta$ the \textbf{page}.
    Note that $\pp P_\theta=B$ for any $\theta$.
    \end{itemize}
Let $X$ be a vector field on $Y$, and $(B,\pi)$ be an open book on $Y$. If $X$ is transverse to each page $P_\theta$ and tangent to $B$, then we say $X$ \textbf{is adapted to} $(B,\pi)$.

\begin{lmm}\label{Open_Book_Lemma}
Let $Y$ be a closed manifold, $B \subset Y$ be a codimension 2 closed submanifold, $X$ be a vector field on $Y$, and $\pi$ be a map from $Y\setminus B$ to $S^1$.
Assume that
\begin{enumerate}[]
    \item $d_p \pi(X) > 0$ for any $p \in Y$,
    \item $B$ has a trivial tubular neighborhood, say $B \times D^2$, such that $\pi(b;r,\theta) = e^{i\theta}$,
    \item $X$ is tangent to $B$.
\end{enumerate}
Then $(B,\pi)$ is an open book decomposition of $Y$, which $X$ is adapted to.
\end{lmm}

\begin{proof}
By (1), $\pi$ is a submersion.
We can remove an open tubular neighborhood $\nu(B)$ of $B$ on which (2) holds, then $\pi|_{Y\setminus\nu(B)}$ is still a submersion.
In particular, (2) guarantees that $\pi|_{\pp\nu(B)}$ is a submersion.
Since $\pi|_{Y\setminus\nu(B)}$ is proper, we can apply the Ehresmann fibration theorem and conclude that $\pi$ defines a fiber bundle. With (2), we can see that $(B,\pi)$ is an open book decomposition on $Y$.

Since $d\pi(X) \neq 0$, $X$ cannot be tangent to the level sets of $\pi$. This means that $X$ is transverse to each page. With (3), we can see that $X$ is adapted to $(B,\pi)$.
\end{proof}

If a vector field $X$ is adapted to $(B,\pi)$, then each page $P_\theta = \pi^{-1}(e^{i\theta})$ can be regarded as a candidate for the global hypersurface of section.
There might exist an orbit of $X$ which does not return to the page in a finite time. Such an orbit should be asymptotic to the boundary as $t$ becomes large.
A discussion about a case of dimension 3 can be found in \cite{Hofer_Wysocki_Zehnder_98}.
In \Cref{Existence_of_GSS_1}, we will show that the case of unbounded return time does not appear in our setting.

\subsection{Symplectic Manifolds and Contact Manifolds}

Let $W$ be a manifold without boundary and $\omega$ be a 2-form on $W$. If $\omega$ is a non-degenerate closed form, we say $\omega$ is a \textbf{symplectic form} and $(W,\omega)$ is a \textbf{symplectic manifold}.
A \textbf{symplectomorphism} is a diffeomorphism between symplectic manifolds which preserves the symplectic form.
Let $W'\subset W$ be a submanifold such that $\omega|_{W'}$ is also a symplectic form on $W'$. We call $W'$ a \textbf{symplectic submanifold}.
A diffeomorphism between manifolds induces a symplectomorphism between cotangent bundles via pullback.
Also, if $N$ is a submanifold of $M$, then $T^*N$ is a symplectic submanifold of $T^*M$.

Let $(W,\omega)$ be a symplectic manifold with a boundary, and assume that $\omega$ has a primitive $\lambda$, i.e.\,\,$\omega=d\lambda$.
A vector field $X$ such that $i_X\omega=\lambda$ is called \textbf{Liouville vector field}.
If a Liouville vector field $X$ exists and defines an outward vector field on the boundary, we call $(W,\lambda)$ a \textbf{Liouville domain}.

Given a smooth function $H:W\to\R$, we can associate a vector field $X_H$ to $H$ via
$$\left(i_{X_H}\right)\omega(-) = \omega(X_H,-) = dH(-).$$
This is well-defined by non-degeneracy of $\omega$.
We say $X_H$ is a \textbf{Hamiltonian vector field}, and in this sense, we call $H$ a \textbf{Hamiltonian function} or simply \textbf{Hamiltonian}.
If a diffeomorphism can be written as a time 1-flow of a Hamiltonian vector field, we say it is a \textbf{Hamiltonian diffeomorphism}.
Note that we can also use time-dependent Hamiltonian $H:W\times\R\to\R$ and get time-dependent Hamiltonian diffeomorphism.
By definition, a Hamiltonian vector field $X_H$ vanishes at a point $p$ if and only if $d_p H=0$.
A Hamiltonian diffeomorphism is a symplectomorphism, as $\left(Fl^{X_H}\right)^*\omega = \mathcal{L}_{X_H}\omega = di_{X_H}\omega =0$ by the Cartan magic formula.
An important observation is that the Hamiltonian vector field of an autonomous Hamiltonian is tangent to the regular level set of the generating Hamiltonian.
In other words, the value of $H$ is preserved under the flow of $X_H$.
More detailed explanations with examples about symplectic geometry and Hamiltonian mechanics can be found in \cite{Arnold_89}, \cite{da_Silva_01}, \cite{Berndt_01}, \cite{Hofer_Zehnder_11} or \cite{McDuff_Salamon_17}.

Let $H,F:W\to\R$ be Hamiltonians. The \textbf{Poisson bracket} is defined by $\{H,F\} : = \omega(X_H,X_F)$.
It's clear that the Poisson bracket is alternating. If $(W,\omega)=(T^*\R^n,\omega_{\mathrm{std}})$, we have the following formula, which can be found for example in the chapter 1 of \cite{McDuff_Salamon_17}
$$\{H,F\} = \sum_j \frac{\partial H}{\partial y_j}\frac{\partial F}{\partial x_j} - \frac{\partial H}{\partial x_j}\frac{\partial F}{\partial x_j}.$$

Let $Y$ be a $(2n+1)$-dimensional manifold, and $\xi$ be a $2n$-dimensional distribution on $Y$.
Then we can locally write $\xi$ as $\ker\alpha$ for some 1-form $\alpha$.
Assume that $\xi$ is coorientable, i.e.\,\,$TW/\xi$ is orientable.
Then we can find a globally defined 1-form $\alpha$.
If $\alpha\wedge(d\alpha)^n$ is a volume form, we say $\alpha$ is a \textbf{contact form}, $\xi=\ker\alpha$ is a \textbf{contact structure}, and $(W,\xi)$ is a \textbf{contact manifold}.
A contact form $\alpha$ defines a unique vector field $R$ such that $\alpha(R) = 1$, $i_R d\alpha =0$.
We call $R$ the \textbf{Reeb vector field}.
Note that the Reeb vector field is always non-vanishing.
A standard example of a contact manifold is the boundary of a Liouville domain $(W,\lambda)$ with contact form $\lambda|_{\pp W}$, and a regular level set of Hamiltonian function $H:(W,d\lambda)\to \R$ with the same contact form.
In the second case, the Hamiltonian vector field $X_H$ restricted to the regular level set of $H$ is the Reeb vector field.
We can consider an open book decomposition $(B,\pi)$ on a contact manifold $(Y,\ker\alpha)$ which the Reeb vector field adapted to.
The close relationship between contact structures on a manifold and an open book decomposition is explored in \cite{Giroux_02}.

\subsection{Geodesic Flow as a Hamiltonian Flow}

Let $(M,g)$ be a complete Riemannian manifold.
For each $x\in M$ and $v\in T_p M$, we have a unique geodesic $\gamma_{x,v}$ with the initial condition $\gamma_{x,v}(0)=x$, $\dot{\gamma}_{x,v}(0)=v$.
The \textbf{geodesic flow} is a 1-parameter family of diffeomorphisms on $TM$ defined by $\Phi_t(x,v)=(\gamma_{x,v}(t),\dot{\gamma}_{x,v}(t))$.
By differentiating $\Phi_t$ by $t$, we get the \textbf{geodesic vector field} on $TM$ which generates the geodesic flow.

The metric $g$ on $TM$ induces the dual metric $g^*$ on $T^*M$ by natural pairing, and we can also define the (co-)geodesic flow and (co-)geodesic vector field on $T^*M$.

\begin{prop}
The geodesic vector field on $T^*M$ is a Hamiltonian vector field with Hamiltonian
    $$
    H(x,y) := \frac{1}{2}\left(||y||_{g^*}^2-1\right).
    $$
where $T^*M$ is equipped with a canonical symplectic form $\omega=\sum dx_i\w dy_i$
\end{prop}

\begin{proof}
    See the proof of Theorem 2.3.1 in \cite{Frauenfelder_van_Koert_18} or Theorem 1.5.2 in \cite{Geiges_08}.
\end{proof}

The regular level set $H^{-1}(0)$ is a unit cotangent bundle $ST^*M$, which naturally is a contact manifold, whose Reeb vector field is Hamiltonian vector field $X_H$ restricted to $ST^*M$.

For notational convenience, we write a point $(x_1,\ldots,x_n)$ in $\R^n$ by $\vec{x}$ so that $x=(x_0,\vec{x})$ is a point in $\R^{n+1}$.
For a function $f:\R^{n+1}\to\R$, we denote the partial derivative $\partial_{x_i}f$ by $f_i$, and the gradient vector field $(f_0,\ldots,f_n)$ by $\nabla f$.
We will also write the Hessian matrix of $f$ which is computed in $\R^{n+1}$ by $\text{Hess}(f) = (\partial_{x_i}\partial_{x_j}f)_{ij} = (f_{ij})$.
We use coordinate $(y_0,\cdots,y_n)$ for the cotangent fiber $T^*_x \R^{n+1}$, and also for the cotangent fiber of a hypersurface contained in $\R^{n+1}$.

Let $f:\R^{n+1}\to\R$ be a smooth function, and $0$ be a regular value.
The level set $M=f^{-1}(0)$ is an $n$-dimensional Riemannian manifold, whose metric is inherited from $\R^{n+1}$.
We can embed $T^* M$ into $T^* \R^{n+1}$ by
    $$T^*M = \left\{(x,y)\in T^*\R^{n+1} : f(x)=0,\,\,y\cdot\nabla f=0\right\}.$$
Here, $\cdot$ is a standard inner product on $\R^{n+1}$, and we identified $T^*M$ to $TM$ by the metric on $M$.
Let $\tilde{H}=\frac{1}{2}\left(||y||^2-1\right)$ on $T^*\R^{n+1}$, and $H=\tilde{H}|_{T^*W}$. Define $\tilde{f},g:T^*\R^{n+1}\to\R$ by
    $$\tilde{f}(x,y)=f(x),\quad g(x,y)=y\cdot\nabla f$$
so that $T^*M$ is an intersection $\tilde{f}^{-1}(0)\cap g^{-1}(0)$.

\begin{prop}\label{Hamiltonian_Vector_Field_on_Submanifold}
    Let $W$ be a symplectic manifold, and $\tilde{H}$ be a Hamiltonian on $W$.
    Consider smooth functions $f,g:W\to\R$ with $c_1,c_2$ as the regular values of $f,g$ such that $V=f^{-1}(c_1)\cap g^{-1}(c_2)\subset W$ is a symplectic submanifold of codimension 2.
    Let $H=\tilde{H}|_V$.
    Then the Hamiltonian vector field $X_H$ is given by
    $$
    X_H = X_{\tilde{H}} - \frac{\{g,\tilde{H}\}}{\{g,f\}}X_f
    -\frac{\{f,\tilde{H}\}}{\{f,g\}}X_g.
    $$
\end{prop}
\begin{proof}
    We can write $X_H$ as $X_H = X_{\tilde{H}} + aX_f + bX_g$ for some functions $a,b$.
    Since $X_H$ is defined on the level set of $f$ and $g$, we must have $X_H(f)=0=X_H(g)$.
    We also have $X_f(f)=0=X_g(g)$ from the definition.
    It follows that 
    $$
    \begin{aligned}
        0 &= X_{\tilde{H}}(f) + bX_g(f) = \{f,\tilde{H}\} + b\{f,g\},\\
        0 &= X_{\tilde{H}}(g) + aX_f(g) = \{g,\tilde{H}\} + a\{g,f\}.
    \end{aligned}
    $$
    Putting $a,b$ into the first equation yields the result.
\end{proof}

\begin{remark}\rm
    The formula of \Cref{Hamiltonian_Vector_Field_on_Submanifold} can be directly generalized to the case of any $2k$-functions, say $V = f_1^{-1}(c_1)\cap\cdots\cap f_{2k}^{-1}(c_{2k})$.
    The coefficients of $X_{f_i}$'s in that case will be a solution of a linear equation consisting of Poisson brackets.
\end{remark}

From a straightforward computation, we have
$$
X_{\tilde{H}} = \sum_j y_j \frac{\partial}{\partial x_j}, \quad 
X_{\tilde{f}} = -\sum_j \frac{\partial f}{\partial x_j} \frac{\partial}{\partial y_j}.
$$
Using the formula for the Poisson bracket, we have
    $$
    \begin{aligned}
    \{\tilde{f},g\} &=
    \sum \frac{\pp f}{\pp x_j}\frac{\pp g}{\pp y_j} = \sum \left(\frac{\pp f}{\pp x_j}\right)^2 = ||\d f||^2,\\
    \{\tilde{f},\tilde{H}\} &
    = \sum \frac{\pp f}{\pp x_j}\frac{\pp \tilde{H}}{\pp y_j} = \sum\frac{\pp f}{\pp x_j}y_j = y \cdot \d f =0,\\ 
    \{g,\tilde{H}\} &
    = \sum \frac{\pp g}{\pp x_j}\frac{\pp \tilde{H}}{\pp y_j}
    = \sum_j \frac{\pp}{\pp x_j}
    \left(\sum_i y_i\frac{\pp f}{\pp x_i}
    \right) y_j
    = \sum _{i,j} \frac{\pp^2 f}{\pp x_i\pp x_j}y_iy_j = \Hess(f)(y,y).
    \end{aligned}
    $$
To sum up, we have the following.
\begin{thm}\label{Geodesic_Vector_Field_Formula}
    Let $f:\R^{n+1}\to\R$ be a smooth function, and $c$ be a regular value of $f$.
    Let $M=f^{-1}(c)$.
    Then, the geodesic vector field on $T^*M$ is given by
    $$
    X_H = \sum_j y_j \frac{\partial}{\partial x_j} - \frac{\Hess(f)_x(y,y)}{\|\nabla f(x)\|^2} \sum_j f_j(x) \frac{\partial}{\partial y_j}.
    $$
\end{thm}
Alternatively, this formula can be derived by orthogonal projection.

\subsection{Sectional Curvature of Hypersurfaces}

An important quantity associated to the Riemannian manifold $(M,g)$ is the \textbf{sectional curvature} $K_M(-,-)$, which is defined for a point $x\in M$ and two vectors $v,w$ in $T_p M$. This can be regarded as the Gauss curvature of a 2-dimensional subspace `spanned' by $v,w$.
For the precise definition and properties of sectional curvature, we refer to \cite{Milnor_63}, \cite{Kobayashi_Nomizu_63}, \cite{Spivak_79_II}, or \cite{do_Carmo_92}.
Let $N\subset M$ be an oriented hypersurface, and $\nu$ be the unit normal vector field on $N$.

\begin{prop}\label{Sectional_Curvature_of_Hypersurface_1}\rm{[\cite{do_Carmo_92}, Theorem 6.2.5]}
    Let $S(v) = \nabla_v \nu$. Then the following formula holds.
    $$
    K_N(v,w) = K_M(v,w) + \frac{\langle S(v), v \rangle \langle S(w), w \rangle - \langle S(v), w \rangle^2}{\langle v, v \rangle \langle w, w \rangle - \langle v, w \rangle^2}
    $$
\end{prop}

The formula becomes simpler in the case of a hypersurface in the Euclidean space, whose sectional curvature vanishes.
\begin{cor}\label{Sectional_Curvature_of_Hypersurface_2}\rm{[\cite{do_Carmo_92}, Theorem 6.2.5]}
    Let $f: \R^{n+1} \to \R$ be a function, and $M = f^{-1}(c)$ be a regular hypersurface. Let $v,w$ be orthonormal vector fields on $TM$. Then we have the following.
    $$
    K_M(v,w) = \langle S(v), v \rangle \langle S(w), w \rangle - \langle S(v), w \rangle^2.
    $$
\end{cor}

\section{A Global Hypersurface of Section of the Geodesic Flow on a Convex Hypersurface}

\subsection{Setting}

We consider the case $f: \R^{n+1} \to \R$ is a smooth function and $0$ is a regular value. Let $M$ be the regular level set  $f^{-1}(0)$. We assume the following conditions for $f$.

\begin{enumerate}[label=\textrm{(A\arabic*)}]
    \item For any point $(x_0, \vec{x})$ in $M$ such that $|x_0|$ is small enough, $f(x_0, \vec{x}) = f(-x_0, \vec{x})$. \label{A1}
    \item The Hessian of $f$ is positive definite.
    Explicitly, for any $(x, y) \in T^*M$, $y \neq 0$, $\Hess(f)_x(y, y) > 0$.\label{A2}
    \end{enumerate}
Here is a classical theorem which explains the role of \ref{A1}.
For the proof, see \cite{Kobayashi_Nomizu_69} Chapter VII.8.

\begin{lmm}\label{Isometric_involution_totally_geodesic_submanifold}
    Let $(M,g)$ be a Riemannian manifold and $N$ be a closed submanifold.
    Assume that there exist a tubular neighborhood $\nu(N)$ of $N$ and an isometric involution $i: \nu(N) \to \nu(N)$, i.e., $i$ is an isometry, $i\neq\mathrm{Id}$ and $i^2 = \mathrm{Id}$.
    In addition, assume that $N$ is the fixed point locus of $i$.
    Then $N$ is a totally geodesic submanifold.
    In other words, the geodesic vector field of $T^*M$ restricted to $T^*N$ is tangent to $T^*N$.
\end{lmm}
Let $N$ be the codimension $1$ submanifold
$N = M \cap \{(x,y): x_0 = 0\}.$
Then it is clear that $N$ is a fixed point set of a locally defined isometric involution $i_M(x_0, \vec{x}) = (-x_0, \vec{x})$, and it follows that $N$ is a totally geodesic submanifold from \Cref{Isometric_involution_totally_geodesic_submanifold}.
Let $Y = ST^*M$ and $B = ST^*N$. Then $B$ can be written as
$$B = \{(x,y) \in ST^*M : x_0 = 0, y_0 = 0\}.$$
We note the triviality of the normal bundle of $B$ in $Y$.

\begin{lmm}\label{Trivial_Normal_Bundle}
    Let $M$ be a Riemannian manifold diffeomorphic to an $n$-sphere, and $N \subset M$ be a codimension $1$ closed submanifold.
    Then $ST^*N$ has a trivial normal bundle in $ST^*M$.
\end{lmm}

\begin{proof}
    We'll first show that the normal bundle $\nu(N)$ of $N$ is trivial, which is equivalent to the orientability of $N$.
    Suppose not, then for any section $s$ of $\nu(N)$ which meets zero section transversely, there exists a loop $\gamma$ such that $s|_{\gamma}$ has odd number of zeros.
    It follows that the intersection form $[N]\cdot[\gamma]$ in $\Z_2$-coordinate is not $0$, which is a contradiction because $H_{n-1}(S^n;\Z_2)=H_1(S^n;\Z_2)=0$.

    Now, let $\nu_M(N) \simeq N \times (-\epsilon, \epsilon) \subset M$. For $x \in N$, we have $T_x M = \mathbb{R} \oplus T_x N$, implying
    \[
    ST_x M = \left\{(t,v) \in \mathbb{R} \oplus T_x N : t^2 + \lVert v \rVert^2 = 1\right\},
    \]
    where $\lVert \cdot \rVert$ is the metric inherited from $M$. Thus, the normal fiber at $p = (x,v) \in ST^*N$ is
    \[
    \nu_{ST^*M}(ST^*N)_p = \nu_{ST_x^*M}(ST_x^*N)_v \oplus \nu_M(N)_x \simeq \nu_{S^{n-1}}(S^n)_v \oplus \nu_M(N)_x,
    \]
    where $S^{n-1}$ is embedded in $S^n$ along the equator.
\end{proof}

Condition \ref{A2} implies that $f$ is a convex function, and $M$ bounds a compact convex domain.
It follows that $M$ is diffeomorphic to $S^n$, where the diffeomorphism can be explicitly written as $x\mapsto x/||x||$, assuming that $0$ is in the interior of the bounding region of $M$.
Since $M$ is convex and $0$ is in the interior of the domain $M$ bounds, this is well-defined.

\subsection{Existence of a Global Hypersurface of Section}

With the setting of the previous section, we define a map $\pi: Y \setminus B \to S^1 \subset \mathbb{C}$ by
$$\pi(x,y) = \frac{x_0 + iy_0}{|x_0 + iy_0|}.$$
The \textbf{angular form} is defined by
$$\Theta = i\cdot d\log\pi = \frac{y_0 dx_0 - x_0 dy_0}{x_0^2 + y_0^2} = \frac{\theta}{x_0^2 + y_0^2}.$$
Put $X_H$ which was computed in \Cref{Geodesic_Vector_Field_Formula} into $\theta$, we have the following.
$$\theta(X_H) = y_0^2 + x_0^2 \frac{\Hess(f)_x(y,y)}{\|\nabla f(x)\|^2} \frac{f_0(x)}{x_0}.$$
We define a function $A = A(x,y)$ by
$$A(x,y) = \frac{\Hess(f)_x(y,y)}{\|\nabla f(x)\|^2}\frac{f_0(x)}{x_0}$$
so that $\theta(X_H) = A(x,y)x_0^2 + y_0^2$.

\begin{prop}\label{Angular_Form_Bound}
Under the assumptions \ref{A1} and \ref{A2}, there exists $\varepsilon > 0$, which only depends on the function $f$, such that $A(x,y) > \varepsilon$ for any $(x,y) \in Y \setminus B$.
\end{prop}

\begin{proof}
    Since $0$ is a regular value of $f$ and $Y$ is compact, there exists some $C > 0$ such that $0<||\d f||^2 \leq C$ on $Y$.
    Also, by the compactness of $Y$, condition \ref{A2} implies that there exists some $\delta > 0$ such that $\mathrm{Hess}(f)_x(y, y) > \delta$ for any $(x,y)\in Y$.

    By condition \ref{A1}, we have $f_0|_{x_0=0} = 0$.
    Thus, we can perform a Taylor expansion with respect to $x_0$:
    \begin{align*}
        f(x_0, \vec{x}) & = f(0, \vec{x}) + \frac{1}{2} f_{00}(0, \vec{x}) x_0^2 + O(x_0^3), \\
        \frac{f_0(x_0, \vec{x})}{x_0} & = f_{00}(0, \vec{x}) + O(x_0).
    \end{align*}
    Since $\Hess(f)_x(y,y)>\delta$, we have $f_{00}(x) = \Hess(f)_x((y_0,0),(y_0,0)) > \delta$ for any $x$. We can choose a small $\eta > 0$ such that for $|x_0| < \eta$, $f_0/x_0 > \delta/2$. Consequently, we have
    $A(x, y) > \delta^2/2C$ for $|x_0|<\eta$.

    For $|x_0| \geq \eta$, we only need to bound $f_0/x_0$.
    However, $f_0(0, \vec{x}) = 0$ and $f_{00} > \delta$ implies that $f_0(x_0, \vec{x}) > 0$ if $x_0 > 0$ and $f_0(x_0, \vec{x}) < 0$ if $x_0 < 0$.
    By the compactness of $Y \cap \{|x_0| \geq \eta\}$, there exists some $\delta_1 > 0$ such that $f_0/x_0 > \delta_1$. This implies $A(x, y) > \delta\delta_1/C$ for $|x_0|\geq\eta$.
    Taking $\varepsilon = \min(\delta^2/2C, \delta\delta_1/C)$, we obtain the desired lower bound.
\end{proof}

\begin{thm}\label{Open_Book_Existence_1}
    Under the assumptions \ref{A1} and \ref{A2}, $\pi: Y \setminus B \to S^1$ defines an open book decomposition, which the geodesic vector field is adapted to.
\end{thm}

\begin{proof}
    We will apply \Cref{Open_Book_Lemma} to this situation. From \Cref{Angular_Form_Bound}, it is evident that there exists $\varepsilon > 0$ such that
    \[
    \Theta(X_H) = \frac{A(x, y)x_0^2 + y_0^2}{x_0^2 + y_0^2} > \varepsilon
    \]
    for any $(x, y) \in Y\setminus B$. Since $\Theta = i \cdot d\log \pi$, the first condition on $\pi$ and $X_H$ in \Cref{Open_Book_Lemma} is satisfied.

    Now consider the trivial tubular neighborhood of $B$, the existence of which is guaranteed by \Cref{Trivial_Normal_Bundle}, denoted as
    \[
    \begin{aligned}
    \nu(B) &\simeq B\times D^2 \\
    (x, y) &\mapsto (\vec{x}, \vec{y}; x_0, y_0)
    \end{aligned}
    \]
    where $x_0, y_0$ are small enough. Note that $\pi(b, r, \theta) = e^{i\theta}$.

    Lastly, since $N$ is a totally geodesic submanifold by \ref{A1}, the geodesic vector field $X_H$ is tangent to $ST^*N = B$.
    Therefore, we can apply \Cref{Open_Book_Lemma} and obtain the desired result.
\end{proof}

\begin{thm}\label{Existence_of_GSS_1}
    Under the assumptions \ref{A1} and \ref{A2}, the geodesic flow on $Y = ST^*M$ admits a global hypersurface of section, which is given by
    $$P = \{(x,y) \in Y : x_0 = 0, y_0 \geq 0\}.$$
\end{thm}

\begin{proof}
    The hypersurface $P$ corresponds to a page $\pi^{-1}(i)$ of the open book constructed in \Cref{Open_Book_Existence_1}.
    To conclude, we need to demonstrate that the return time is bounded.
    It suffices to show that there exists $t\in\R_{>0}$ such that $\pi(Fl^{X_H}_t(x, y)) = \pi(x, y)$.
    From \Cref{Angular_Form_Bound}, we observe the existence of $\varepsilon > 0$ such that $\Theta(X_H) > \varepsilon$.
    Thus, we have
    \[
    \int_0^{2\pi/\varepsilon} i\cdot d\log\pi(X_H) > \int_0^{2\pi/\varepsilon} \varepsilon dt > 2\pi.
    \]
    By the intermediate value theorem, we can conclude that there exists a positive number $\tau<2\pi/\varepsilon$ such that $\int_0^\tau i\cdot d\log\pi(X_H)=2\pi$, which means that $\pi(Fl^{X_H}_\tau(x,y))=\pi(x,y).$
    It means that $\tau$ is a bounded positive finite return time for $(x,y)$.
    The boundedness of the negative return time can be demonstrated similarly.
\end{proof}

\begin{cor}
    Under the assumptions \ref{A1} and \ref{A2}, the geodesic flow on $M$ admits an $S^1$-family of global hypersurfaces of section.
\end{cor}

\begin{proof}
    We have an open book decomposition with bounded return time, so we can choose any page of the open book, let's say $\pi^{-1}(e^{i\theta}) = P_\theta$, as a global hypersurface of section.
\end{proof}

Note that if we take $n=2$, then $M$ is 2-sphere so $P=P_{\pi/2}$ in \Cref{Existence_of_GSS_1} is the Birkhoff annulus.
In this sense, our construction can be regarded as a generalization of the Birkhoff annulus.

\subsection{Relation with Sectional Curvature}

We can apply \Cref{Sectional_Curvature_of_Hypersurface_2} on our hypersurface $M$.
The unit normal vector $\nu$ is $\frac{\d f}{||\d f||}$, and we have
$$
S\left(\frac{\pp}{\pp x_i}\right) = \frac{\pp}{\pp x_i}\frac{\d f}{||\d f||} = \sum_j \left( \frac{f_{ij}}{||\d f||} - \frac{\sum_k f_k f_{ki} f_j}{||\d f||^3}\right)\frac{\pp}{\pp x_j}.
$$
It follows that
$$
S(v,w) = \frac{\sum_{i,j} v_i f_{ij} w_i}{||\d f||} - \frac{\sum_{i,j,k} v_i f_{ik} f_k f_j w_j}{||\d f||^3} = \frac{\Hess(f)(v,w)}{||\d f||}
$$
for $v,w\in TM$. The second term vanishes because $w\in TM$ implies that $w\cdot \d f=\sum_j w_j f_j = 0$.

Put this into \Cref{Sectional_Curvature_of_Hypersurface_2}, and we get the following.

\begin{prop}\label{Sectional_Curvature_and_Hessian}
    Let $v,w$ be orthogonal unit tangent vectors of a point $x \in M$. Then,
    $$
    K_M(v,w) = \frac{\Hess(f)_x(v,v)\Hess(f)_x(w,w) - \Hess(f)_x(v,w)^2}{||\d f(x)||^2}.
    $$
\end{prop}

Hence, the sign of $K_M$ depends on the sign of the ($2\times2$)-minors of $\Hess(f)$. Since the formula holds for any $v$ and $w$, we have the following.

\begin{cor}\label{Curvature_and_Hessian}
The sectional curvature $K_M(\sigma)_x$ is positive for any point $x\in M$ and plane $\sigma\subset T_x M$ if and only if $\Hess(f)_x$ is either positive definite or negative definite.
\end{cor}
\begin{proof}
Since $\Hess(f)_x$ is symmetric bilinear form, there exists a basis $\mathfrak{B}$ of $T_x M$ that diagonalizes $\Hess(f)_x$, and it continuously varies with $x$. Let
$$[\Hess(f)_x]_\mathfrak{B} = \mathrm{diag}(\lambda_1, \cdots, \lambda_n).$$
Since $K_M > 0$, we must have $\lambda_i \lambda_j > 0$ for any $i,j$. It means that $\lambda_i$ all have the same sign, which implies that $\Hess(f)_x$ is either positive definite or negative definite. Since the process depends on $x$ continuously, the sign of $\lambda_i$ cannot change since $\lambda_i \neq 0$. The converse follows from the formula in \Cref{Sectional_Curvature_and_Hessian}.
\end{proof}

The following corollary can also be found in various literature, for example \cite{Sacksteder_60}.

\begin{cor}\label{Convex_Hypersurface_is_Sphere}
    Let $f:\R^{n+1}\to\R$ and $M = f^{-1}(0)\subset \R^{n+1}$ be a regular hypersurface. Assume that $M$ has a positive sectional curvature. Then $M$ is diffeomorphic to the $n$-sphere.
\end{cor}
\begin{proof}
By \Cref{Curvature_and_Hessian}, $\Hess(f)$ is either positive definite or negative definite. If $\Hess(f)$ is positive definite, $f^{-1}(-\infty,0]$ is convex and we can get the result. If $\Hess(f)$ is negative definite, then for $\bar{f}=-f$, $M=\bar{f}^{-1}(0)$, and $\Hess(\bar{f})$ is positive definite, and we can use the same argument.
\end{proof}

Now we can formulate \Cref{Existence_of_GSS_1} in another form.

\begin{thm}\label{Open_Book_Existence_2}
    Let $M = f^{-1}(0)\subset\R^{n+1}$ be a hypersurface satisfying \ref{A1}, and $M$ has a positive sectional curvature. Then there exists an open book decomposition $(\pi, B)$ of $ST^*M$ to which the geodesic vector field $X_H$ is adapted, as given in \Cref{Open_Book_Existence_1}. Moreover, there exists a global hypersurface of section $P$, as given in \Cref{Existence_of_GSS_1}.
\end{thm}
\begin{proof}
Since $K_M > 0$, \Cref{Curvature_and_Hessian} implies that $f$ is either positive definite or negative definite. If $f$ is positive definite, we're done. If $f$ is negative definite, define $\tilde{f}=-f$. Then $\tilde{f}$ is positive definite, $\tilde{f}^{-1}(0) = f^{-1}(0)$, and the condition \ref{A1} still holds. Moreover, the formula in \Cref{Geodesic_Vector_Field_Formula} for $f$ and $\tilde{f}$ are the same. Thus we get the result.
\end{proof}

\subsection{Topology of the Global Hypersurface of Section}

Using the same notations from the previous sections, we will explore the topology of $N$ and $P$.

\begin{lmm}\label{Topology_of_Hypersurface}
Under the assumptions \ref{A1} and \ref{A2}, $N$ is diffeomorphic to $S^{n-1}$.
\end{lmm}

\begin{proof}
We can express $N$ as
$$
N = \left\{
(0,\vec{x})\in\R^{n+1}:
f(0,\vec{x}) = 0
\right\}
\subset \left\{
(0,\vec{x})\in\R^{n+1}:
\vec{x}\in\R^n
\right\}\simeq\R^n.
$$
The convexity of $f$ is preserved if we restrict $f$ to the subspace $\left\{(0,\vec{x}):\vec{x}\in\R^n\right\}$, so we obtain the result.
\end{proof}

\begin{prop}\label{Topology_of_Page}
The global hypersurface of section $P$ constructed in \Cref{Existence_of_GSS_1} is diffeomorphic to $T_{\leq1}^*S^{n-1}$, which is the subset of $T^*S^{n-1}$ consisting of covectors of length $\leq1$.
The boundary $B=\pp P$ is homeomorphic to $ST^*S^{n-1}$.
\end{prop}

\begin{proof}
Let $p$ be a diffeomorphism from the upper hemisphere $H^n\subset S^n$ to the closed disk $D^n$.
Define a map $\phi: P\to T^*_{\leq1}N$ by $\phi(x,y)=(x,p(y))$, then it's clear that $\phi$ is a diffeomorphism.
With \Cref{Topology_of_Hypersurface}, we can conclude the result.
\end{proof}

\section{Return Map}

\subsection{General Properties}

In this subsection, we investigate the general properties of the global hypersurfaces of sections for Reeb vector fields.
Let $\alpha$ be a contact form on a manifold $Y$ and $R$ be its Reeb vector field.
Assume that there exists an open book decomposition $(B,\pi)$ of $Y$ to which $R$ is adapted.
Let $P=P_\theta=\pi^{-1}(e^{i\theta})$ be a page, and assume that the first return map $\Psi:\mathring{P}\to \mathring{P}$ is well-defined, i.e.\,\,the return time $\tau$ is bounded for each point $p\in P$.

\begin{lmm}\label{Page_is_Symplectic}
The interior of $P$ is a symplectic manifold with symplectic form $d\alpha$.
\end{lmm}

\begin{proof}
Let $\dim Y=2n+1$. Since $Y$ is contact, $\alpha\w d\alpha^n$ never vanishes. Since $\mathring{P}$ is transverse to $R$, we can take a local frame $(X_1,\cdots,X_{2n})$ of $\mathring{P}$ such that $(R,X_1,\cdots,X_{2n})$ is a local frame of $Y$, and
$$
\alpha\w d\alpha^n(R,X_1,\cdots,X_{2n}) =\alpha(R)d\alpha^n(X_1,\cdots,X_{2n})\neq0.
$$
It means that $d\alpha$ is a non-degenerate closed 2-form on $\mathring{P}$, so it's a symplectic form.
\end{proof}

\begin{prop}\label{Return_Maps_is_Symplectic}
For the return map $\Psi:\mathring{P}\to\mathring{P}$, we have
$$
\Psi^*\alpha-\alpha=d\tau.
$$
In particular, $\Psi$ is a symplectomorphism.
\end{prop}

\begin{proof}
This is a generalization of the well-known fact in dimension 3, which can be found, for example, in \cite{Abbondandolo_Bramham_Hryniewicz_Salomao_17} or \cite{Frauenfelder_van_Koert_18}. We have $Fl^{R}_{\tau_p}(p) = \Psi(p)$. Differentiating both sides and plugging in a vector field $X$, we have
$$
dFl^R_{\tau_p}(p)X + (d\tau_p(X))R = d_p\Psi (X).
$$
Since $Fl^R$ preserves $\alpha$, we have that
$$
\Psi^*\alpha(X) =\alpha(d\Psi(X)) =(Fl^R)^*\alpha(X) +\alpha(R)d\tau(X) =\alpha(X)+d\tau(X).
$$
The second statement follows if we take exterior derivative on both sides.
\end{proof}

\begin{rmk}\rm
If we equip $d\alpha$ to the whole $P$, which is a manifold with boundary, then $(P, d\alpha)$ is not a Liouville domain.
More precisely, $d\alpha$ degenerates at the boundary of $P$.
The Reeb vector field $R$ is tangent to $\partial P$ and $i_R d\alpha = 0$, indicating that $d\alpha$ degenerates on the boundary.
But even though $d\alpha$ degenerates at the boundary, $(\partial P, \ker(\alpha|_{\partial P}))$ is a contact manifold.
\end{rmk}

\subsection{Extension to the Boundary}

In the preceding sections, we constructed a global hypersurface of section $P$ for the geodesic flow on a convex hypersurface $M = f^{-1}(0)$ in $\R^{n+1}$ and we obtained a return map $\Psi: \mathring{P} \to \mathring{P}$.
As demonstrated in \Cref{Page_is_Symplectic}, $\mathring{P}$ is a symplectic manifold, with its symplectic form being the restriction of the standard symplectic form on $T^*M$. Furthermore, as indicated in \Cref{Return_Maps_is_Symplectic}, $\Psi$ is a symplectomorphism.

Now we investigate the boundary behavior of the return map in terms of the defining function $f$.
Given that $(\partial_{x_0},\partial_{y_0})$ forms a normal symplectic frame of a contact submanifold $B=ST^*N$, it follows from Lemma 8.1 of \cite{Moreno_van_Koert_22a} the linearized Hamiltonian flow decomposes into two blocks; one corresponding to $TB$ and the other to $\langle \partial_{x_0}, \partial_{y_0} \rangle$.
Let $Z_N=\begin{pmatrix}
    x_0(t)\\y_0(t)
\end{pmatrix}$ be the normal part of a flowline linearized near $B$.
Then we have
    $$\dot{Z}_N = J_NS_NZ_N$$
where $J_N$ is the normal part of the standard complex structure on $TY=TT^*M$ and $S_N$ is a matrix called the \textbf{normal Hessian}.

\begin{prop}\label{Extension_to_Boundary_Lemma}
{\rm \cite{Moreno_van_Koert_22a}}
If the normal Hessian is positive definite, the return map extends to the boundary smoothly, and such an extension is unique.
\end{prop}

\begin{thm}\label{Smooth_Extension_to_Boundary}
Under the assumptions \ref{A1} and \ref{A2}, or the assumption of \Cref{Open_Book_Existence_2}, the return map $\Psi: \mathring{P} \to \mathring{P}$ can be smoothly extended to the boundary.
\end{thm}
\begin{proof}
Instead of working with a normal Hessian, we will use linearized flow.
According to Proposition 8.2 of \cite{Moreno_van_Koert_22a}, we must have
    $$
    \Theta(X_H) = \frac{Z_N^t S_N Z_N}{Z_N^t Z_N} + O(1)
    $$
near $(x_0,y_0)=(0,0)$.
From the previous section, we also have
$$\Theta(X_H) = \frac{A(x,y)x_0^2 + y_0^2}{x_0^2 + y_0^2} =\frac{Z_N^t \mathrm{diag}\left(
    \frac{\Hess(f)_x(y,y)}{||\d f(x)||^2}f_{00}(x),1
    \right) Z_N}{Z_N^t Z_N} + O(x_0).
    $$
It follows that
    $$S_N = \mathrm{diag}\left(
    \frac{\Hess(f)_x(y,y)}{||\d f(x)||^2}f_{00}(x),1
    \right).$$
From the condition \ref{A2}, we observe that $\Hess(f)(y,y)$ and $f_{00}$ are always positive, and it follows that $S_N$ is positive definite.
From \Cref{Extension_to_Boundary_Lemma}, we get the result.
If we impose the positive sectional curvature condition, it follows that $\Hess(f)(y,y)$ and $f_{00}$ have the same sign, which leads to the same result.
\end{proof}

\begin{ex}\rm
\Cref{Open_Book_Existence_1} can be applied to ellipsoids, since they are convex and satisfying the symmetry condition \ref{A1}.
In particular, consider an ellipsoid given by
    $$
    E = \left\{x\in\R^{n+1}:\frac{x_0^2}{a_0^2}+x_1^2+\cdots+x_n^2-1=0\right\},
    $$
where $a_0\in\R_{>0}$.
The global hypersurface of section constructed in \Cref{Open_Book_Existence_1} is the upper-hemisphere bundle on the equator,
    $$
    P = \left\{
    ((0,\vec{x}),(y_0,\vec{y}))\in T^*\R^{n+1}\,:\,
    ||\vec{x}||=1,\,\il\vec{x},\vec{y}\ir=0,\,y_0^2+||\vec{y}||^2=1,\,y_0\geq0
    \right\}.
    $$
The return map $\Psi$ on $P$ can be computed in terms of elliptic integrals
        $$
        \Psi\left((0,\vec{x}), (y_0, \vec{y})\right) =
        \left(\left(0, \vec{x}\cos{G(\|\vec{y}\|)} + \frac{\vec{y}}{\|\vec{y}\|}\sin{G(\|\vec{y}\|)}\right),
        \left(y_0, \vec{y}\cos{G(\|\vec{y}\|)} - \|\vec{y}\|\vec{x}\sin{G(\|\vec{y}\|)}\frac{}{}\right)\right),
    $$
where $G$ is given by
    $$
    G(t) = -\frac{t(1-a_0^2)}{a_0}F\left(2\pi\,\Big|\,\frac{-(1-a_0^2)(1-t^2)}{a_0^2}\right)
    +\frac{t}{a_0}\Pi\left(1-t^2;2\pi\,\Big|\,\frac{-(1-a_0^2)(1-t^2)}{a_0^2}\right).
    $$
where $F$ is the elliptic integral of the first kind and $\Pi$ is the elliptic integral of the third kind.

This can be regarded as a specific case of a \emph{hypersurface of revolution}.
Let $f:\R^{n+1}\to\R$ satisfy \ref{A1} globally, \ref{A2} and the third condition; if $||\vec{x}||=||\vec{x}\,'||$, then $f(x_0,\vec{x})=f(x_0,\vec{x}\,')$.
For convenience, we assume that $f(0,\vec{x})=0$ if and only if $||\vec{x}||=1$.
The hypersurface of revolution is given by a regular level set $M=f^{-1}(0)$.
In this case, we can parametrize the set $M\cap \left\{(x_0,x_1,0,\cdots,0)\right\}$ by $(a(\phi),\cos\phi,0,\cdots,0)$, where $a$ is a function on $\phi$. Note that in the case of ellipsoid, $a(\phi)=a_0\sin\phi$.
We can compute the return map with a help of the \emph{Clairaut integral} \cite{Arnold_89}, which has a same form as in the case of ellipsoid with
    $$
     G(t):= t\int_{0}^{2\pi} \frac{\sqrt{(1-t^2)\sin^2{\sigma} + \{a'(\arcsin{(\sqrt{1-t^2}\sin{\sigma})})\}^2}}{1-(1-t^2)\sin^2{\sigma}} d\sigma.
    $$
It's straightforward to see that the return map $\Psi$ is a Hamiltonian diffeomorphism generated by $H(x,y)=\left(\int G\right)(||\vec{y}||)$.

Sending $a(\phi)$ to 0, the dynamics converges to a billiard on the unit disk $D=\left\{\vec{x}\in\R^n:||\vec{x}||\leq1\right\}$, which is defined on a set
    $$
    Y_0 = \left\{
    (\vec{x},\vec{y})\in T^*\R^n:||\vec{x}||\leq1,\,||\vec{y}||=1,\,\vec{y}\text{ points inward if }||\vec{x}||=1
    \right\}.
    $$
The global hypersurface of section converges to $P_0 = \left\{(\vec{x},\vec{y})\in Y_0:||\vec{x}||=1\right\}$.
We write $\vec{y}=\vec{y}_T+\vec{y}_N$ for $(\vec{x},\vec{y})\in P_0$, where $\vec{y}_T$ is a component tangent to $\pp D$, and and $\vec{y}_N$ is a component normal to $\pp D$.
Then the return map converges to the second iterate of the high-dimensional billiard map
$$
        \Psi\left(\vec{x}, \vec{y}\right) =
        \left( \vec{x}\cos{G_0(\|\vec{y}_T\|)} + \frac{\vec{y}_T}{\|\vec{y}_T\|}\sin{G_0(\|\vec{y}_T\|)},\,
        \vec{y_T}\cos{G_0(\|\vec{y}_T\|)} - \|\vec{y}_T\|\vec{x}\sin{G_0(\|\vec{y_T}\|)}\right),
    $$
where $G_0(t)=4\arccos t$, which agrees with the function $G$ with $a(\phi)=0$.
\end{ex}

\bibliographystyle{amsalpha}
\bibliography{CH}
\end{document}